\documentclass[10pt,a4paper]{amsart}

\usepackage{amssymb}
\usepackage{amsfonts}
\usepackage{amsmath}
\usepackage{amsthm}

\numberwithin{equation}{section}
\newtheorem{theorem}{Theorem}[section]

\newtheorem{proposition}[theorem]{Proposition}
\newtheorem{lemma}[theorem]{Lemma}
\theoremstyle{definition}

\theoremstyle{remark}
\newtheorem{remark}{Remark}[section]
\newtheorem*{acknowledgments}{Acknowledgments}

\newcommand{\dd}{\mathrm{d}}
\newcommand{\norm}[3]{\left\|#1\right\|^{#2}_{#3}}
\newcommand{\duality}[2]{\left\langle #1 \Big| #2 \right\rangle}
\newcommand{\dist}{\mathrm{dist}}
\DeclareMathOperator{\supp}{supp}

\newcommand{\R}{\mathbb{R}}
\newcommand{\C}{\mathbb{C}}

\newcommand{\N}{\mathbb{N}}
\newcommand{\BP}{X_{\zeta}^{1/2}}
\newcommand{\BM}{X_{\zeta}^{-1/2}}
\newcommand{\BHP}{\dot{X}_{\zeta}^{1/2}}
\newcommand{\BHM}{\dot{X}_{\zeta}^{-1/2}}

\newcommand{\y}{\mathbf{y}}

\title[]{Stability of the Calder\'on problem for less regular conductivities}
\author[]{Pedro Caro \and Andoni Garc\'ia \and Juan Manuel Reyes}
\date{May 11, 2012}
\address{Department of Mathematics and Statistics, Helsingin yliopisto / Helsingfors universitet / University of Helsinki, Finland}
\email{pedro.caro@helsinki.fi}
\email{andoni.garcia@helsinki.fi}
\email{juan.reyes@helsinki.fi}

\begin{document}

\begin{abstract} In these notes we prove log-type stability for the Calder\'on problem with conductivities in $ C^{1,\varepsilon}(\overline{\Omega}) $. We follow the lines of a recent work by Haberman and Tataru in which they prove uniqueness for $ C^1(\overline{\Omega}) $.
\end{abstract}

\maketitle

\tableofcontents
\setcounter{tocdepth}{1}

\section{Introduction}
Let $ \Omega $ be a bounded non-empty open subset of $ \R^n $ with $ n \geq 2 $ and let $ \partial \Omega $ denote its boundary. Let $ \gamma $ belong to $ L^\infty(\Omega) $ and assume that $ \gamma (x) \geq \gamma_0 > 0 $ for almost every $ x \in \Omega $. Define the Dirichlet-to-Neumann map $ \Lambda_\gamma : H^{1/2} (\partial \Omega) \longrightarrow H^{-1/2} (\partial \Omega) $ as
\[ \duality{\Lambda_\gamma f}{g} := \int_\Omega \gamma \nabla u \cdot \nabla v \, \dd x \]
for any $ f, g \in H^{1/2} (\partial \Omega) $, where $ u \in H^1(\Omega) $ is the weak solution of the conductivity equation $ \nabla \cdot (\gamma \nabla u) = 0 $ in $ \Omega $ with $ u|_{\partial \Omega} = f $ and $ v \in H^1(\Omega) $ with $ v|_{\partial \Omega} = g $. $ \Lambda_\gamma $ is well-defined, linear and bounded. Let $ \norm{\centerdot}{}{} $ denote the operator norm from $ H^{1/2} (\partial \Omega) $ to $ H^{-1/2} (\partial \Omega) $. In this context, the coefficient $ \gamma $ models an isotropic conductivity on $ \Omega $. In \cite{C}, Calder\'on posed the problem consisting in recovering the conductivity $ \gamma $ from the boundary measurements $ \Lambda_\gamma $. This inverse boundary value problem (IBVP) is widely known as the Calder\'on problem.

In these notes, we discuss the stability issue of this IBVP in the particular case that $ n \geq 3 $, $ \partial \Omega $ is locally described by the graph of a Lipschitz function and the conductivities lie in the functional space $ C^{1, \varepsilon}(\overline{\Omega}) $ for an arbitrarily small $ \varepsilon $. Along our discussion we follow a recent improvement of the classical method introduced by Sylvester and Uhlmann in \cite{SU} and based on the construction of complex geometric optics solutions (CGOs for short). This new improvement is due to Haberman and Tataru (see \cite{HT}) and it has allowed us to prove the following stability estimate for the Calder\'on problem.

\begin{theorem} \sl Let $ \Omega $ be a bounded non-empty open subset of $ \R^n $ with $ n \geq 3 $. Assume $ \partial \Omega $ to be locally described by the graph of a Lipschitz function. Let $ M, \delta $ and $ \varepsilon $ be real constants such that $ M > 1 $, $ 0 < \delta < 1 $ and $ 0 < \varepsilon < 1 $. Then,
\[ \norm{\gamma_1 - \gamma_2}{}{C^{0, \delta}(\overline{\Omega})} \lesssim \left( \log \norm{\Lambda_{\gamma_1} - \Lambda_{\gamma_2}}{-1}{} \right)^{- \varepsilon^2 (1 - \delta) / (3n^2)} \]
for all $ \gamma_1, \gamma_2 \in C^{1,\varepsilon}(\overline{\Omega}) $ such that $ \gamma_j (x) > 1/M $ for all $ x \in \Omega $ and $ \norm{\gamma_j}{}{C^{1,\varepsilon}(\overline{\Omega})} \leq M $.
\label{th:stability}
\end{theorem}

For the sake of completeness, let us recall that, for a suitable $ a $,
\[ \norm{a}{}{C^{0, \delta}(\overline{\Omega})} = \inf \{ C > 0 : |a(x)| \leq C,\, |a(x) - a(y)| \leq C |x - y|^\delta\, \forall x,y \in \Omega \}. \]

Along these notes, the symbol $ \lesssim $ means that there exists a positive constant for which the estimate (for the symbol $ \leq $) holds whenever the right hand side of the estimate is multiplied by that constant. We will refer to this constant as implicit constant. In Theorem \ref{th:stability}, the implicit constant just depends on $ n, \Omega, M, \delta $ and $ \varepsilon $.

Since \cite{C}, many papers have addressed the questions of uniqueness and stability related to the Calder\'on problem. Nowadays, the picture of this problem seems to have two faces, one considering the problem in dimension $ n = 2 $ and other for dimension $ n \geq 3 $. In dimension $ n = 2 $ the Calder\'on problem was completely solved in \cite{AsP} by Astala and P\"aiv\"arinta and sharp stability results in this framework where given in \cite{ClFR} by Clop, Faraco and Ruiz and in \cite{FRo} by Faraco and Rogers. Some previous results are \cite{N95}, \cite{BU}, \cite{BarBaR} and \cite{BaFR}. In dimension $ n \geq 3 $ there are still open questions about the optimal smoothness for uniqueness and stability. The best known positive result for uniqueness is due to Haberman and Tataru. They proved in \cite{HT} uniqueness for the Calder\'on problem for continuously differentiable conductivities. Some other previous results are the following: the foundational \cite{SU} by Sylvester and Uhlmann for smooth conductivities and \cite{B} by Brown for conductivities in $ \cup_{\varepsilon > 0} C^{1,1/2+\varepsilon} (\overline{\Omega}) $. Other references to be mentioned are \cite{BTo} by Brown and Torres (uniqueness for conductivities with $ 3/2 $ derivatives in $ L^p $, $ p > 2n $), \cite{PPaU} by P\"aiv\"arinta, Panchenko and Uhlamnn (uniqueness for conductivities in $ W^{3/2,\infty} (\Omega) $) and \cite{N} where Nachman provides an algorithm to reconstruct a conductivity from its Dirichlet-to-Neumann map. As far as we know, the best stability result is due to Heck and it was stated in \cite{He} for conductivities in $ H^{n/2 + \varepsilon}(\Omega) \cap C^{1,1/2+\varepsilon} (\overline{\Omega}) $ and with $ \Omega $ having a smooth boundary. Heck's paper follows the lines of \cite{PPaU} and \cite{A} by Alessandrini (as far as we know, this is the first paper proving internal stability for the Calder\'on problem). It seems that the only regularity assumption in \cite{He} imposed by the method of uniqueness in \cite{PPaU} is to have conductivities in $ C^{1,1/2+\varepsilon} (\overline{\Omega}) $. The extra regularity assumption (namely $ H^{n/2 + \varepsilon}(\Omega) $) seems to be used to control the $ L^\infty $-norm of the conductivities by the boundary data. This idea goes back to \cite{A}.

Haberman and Tataru's ideas allowed us to extend Theorem 1 in \cite{BarBaR} to dimension $ n \geq 3 $ and to improve Heck's result relaxing the smoothness of the coefficients and the smoothness of the boundary of the domain. Our argument also allows us to control the norm of the conductivities in $ C^{0, \delta}(\overline{\Omega}) $ by the boundary data without assuming extra regularity of the coefficients, just by paying with a power on the right hand side of the stability estimate. We can get rid of the unpleasant assumption in \cite{He} (conductivities have to belong to $ H^{n/2 + \varepsilon}(\Omega) $) using interpolation in Lebesgue spaces and Morrey's embedding. In \cite{Ca} and \cite{Ca11} the first author proved stability estimates for an IBVP arising in electromagnetism. In these estimates the $ H^1 $-norms of the electromagnetic coefficients were bounded by the boundary data. The same argument of interpolation in Lebesgue spaces and Morrey's embedding provides now, under the same conditions as in \cite{Ca} and \cite{Ca11}, stability estimates controlling the $ C^{0,\delta} $-norms of the coefficients.

As we already mentioned, our approach to study the stability of this IBVP uses the CGOs constructed by Haberman and Tataru in \cite{HT}. It seems that the main idea in this paper is to prove certain decay properties for the remainder of the CGOs for less regular conductivities in certain Bourgain spaces. It was pointed out in \cite{B} the lack of decay for the remainder whenever the conductivity was less regular than $ C^{1,1/2+\varepsilon} (\overline{\Omega}) $. This seems to be the case even using this type of Bourgain spaces (see Section \ref{sec:CGOreview}). The breakthrough in \cite{HT} is to prove the remainder properties of the CGOs making average over parameters associated to the introduced Bourgain spaces. Thus, they are able to prove decay \emph{in average}.

The structure of the paper is as follows: In Section \ref{sec:BOUNDARY-INTERIOR}, we prove an estimate relating the internal electric properties with the boundary measurements. This estimate might be different to the usual ones since the term containing internal electric properties is written in the whole space and not only in $ \Omega $. We will take advantage of this at the end when using the Fourier transform. In order to write the term containing internal electric properties in the whole space, we have to perform appropriate extensions of the coefficients. These extensions are also carried out in the second section. In Section \ref{sec:CGOreview}, we review the construction of CGOs given by Haberman and Tataru and we set up the precise properties we need to prove Theorem \ref{th:stability}. In Section \ref{sec:stability}, we prove Theorem \ref{th:stability} using the estimate from Section \ref{sec:BOUNDARY-INTERIOR} and the solutions reviewed in Section \ref{sec:CGOreview}. Our proof follows the general lines of \cite{A} but it also requires the stability on the boundary proven by Alessandrini in \cite{A90}. The key ingredient in our proof is the use of the solutions constructed by Haberman and Tataru. However, our way to proceed is slightly different to the one followed by them to prove uniqueness. In \cite{HT}, the authors deduced from the decay in average that, for any Fourier frequency, there exists a sequence of solutions with the \emph{good} remainder properties and they use these sequences of solutions. From the point of view of stability this approach does not seem to be very convenient, so instead of doing so, we use directly the decay in average to prove Theorem \ref{th:stability}. Finally, in Section \ref{sec:finalDISCUSSION} we discuss possible improvements of Theorem \ref{th:stability} following the lines of our argument. As a consequence of this discussion we motivate two naive questions.

\begin{acknowledgments} We thank Juan Antonio Barcel\'o for encouraging us to study this problem and Alberto Ruiz for his valuable comments. We also want to thank the anonymous referee for his or her comments which have contributed to improve this manuscript. The authors are supported by the projects ERC-2010 Advanced Grant, 267700 - InvProb and Academy of Finland (Decision number 250215, the Centre of Excellence in Inverse Problems). PC and JMR also belong to the project MTM 2011-02568 Ministerio de Ciencia y Tecnolog\'ia de Espa\~na. AG belongs to the project MTM2007-62186 Ministerio de Ciencia y Tecnolog\'ia de Espa\~na. \end{acknowledgments}

\section{From the boundary to the interior}\label{sec:BOUNDARY-INTERIOR}
In this section we prove an estimate relating the internal electric properties with the boundary measurements. In order to prove this estimate, we will perform appropriate extensions of the coefficients. In \cite{BaFR}, Barcel\'o, Faraco and Ruiz carried out an argument closely related to the one presented in this section. The main difference lies in the smoothness of the functions to be extended. Conductivities in \cite{BaFR} belong to $ C^{0,\varepsilon} (\overline{\Omega}) $ while here they belong to $ C^{1,\varepsilon} (\overline{\Omega}) $. This makes our argument more technical.

Let $ R $ be a positive constant. Along these notes, $ B $ will always denote the open ball of radius $ R $ given by $ B := \{ x \in \R^n : |x| < R \} $. For us, $ C^{0, \delta}(\R^n) $ with $ \delta \in (0, 1] $ denotes the Banach space of H\"older (Lipschitz if $ \delta = 1 $) continuous functions in $ \R^n $. More precisely,
\[ C^{0,\delta}(\R^n) = \{ a : \R^n \longrightarrow \C \, :\, \exists L > 0 \, |a(x)| \leq L, |a(x) - a(y)| \leq L |x - y|^\delta\, \forall x, y \in \R^n \}, \]
whose norm is the smallest $ L $ in the definition. This norm will be denoted by $ \norm{\centerdot}{}{C^{0, \delta}(\R^n)} $. Recall that $ C^{0,1}(\R^n) $ is equivalent to $ W^{1, \infty} (\R^n) $, the space of measurable functions (modulo those vanishing almost everywhere) such that themselves and their first weak partial derivatives are essentially bounded in $ \R^n $. We are now prepared to prove the following lemma:
\begin{lemma} \sl
Consider a real constant $ \gamma_0 $ in the interval $ (0, 1] $ and $ R $ such that $ \overline{\Omega} \subset B $. Let $ \gamma_1 $ and $ \gamma_2 $ be two given functions belonging to $ C^{0,1}(\R^n) $ satisfying $ \supp (\gamma_j - 1) \subset \overline{B} $ and $ \gamma_j (x) \geq \gamma_0 $ for all $ x \in \R^n $ and $ j \in \{ 1, 2 \} $. Then, for any $ u_j \in H^1_\mathrm{loc}(\R^n) $ weak solution of $ \nabla \cdot (\gamma_j \nabla u_j) = 0 $ in $ \R^n $, one has
\begin{gather}
\left| \int_{\R^n} \nabla \gamma_2^{1/2} \cdot \nabla (\gamma_2^{-1/2} v_1 v_2) \, \dd x - \int_{\R^n} \nabla \gamma_1^{1/2} \cdot \nabla (\gamma_1^{-1/2} v_1 v_2) \, \dd x \right| \leq \nonumber \\
\leq \left( \norm{\Lambda_{\gamma_1|_\Omega} - \Lambda_{\gamma_2|_\Omega}}{}{} + \norm{\gamma_1 - \gamma_2}{}{L^\infty(B \setminus \Omega)} \right) \norm{u_1}{}{H^1(B)} \norm{u_2}{}{H^1(B)}. \label{es:Integralbound}
\end{gather}
Here $ v_j \in H^1_\mathrm{loc}(\R^n) $ denotes $ v_j = \gamma_j^{1/2} u_j $. \label{lem:integralBOUND}
\end{lemma}

In this lemma as in the remainder of these notes, $ \Omega $ is fixed and it satisfies the assumptions of Theorem \ref{th:stability}.

Recall that $ H^1_\mathrm{loc}(\R^n) $ stands for the space of locally integrable functions (modulo those vanishing almost everywhere) such that their restriction and the restriction of their gradient to any compact subset of $ \R^n $ are square-integrable.

\begin{proof} Firstly note that $ \gamma_1 (x) = \gamma_ 2 (x) $ for any $ x \in \R^n \setminus B $. Secondly, if $ j,k \in \{ 1, 2 \} $ and $ k \neq j $, one has
\begin{gather*}
\duality{\Lambda_{\gamma_j|_{\Omega}}(u_j|_{\partial\Omega})}{u_k|_{\partial\Omega}} = \int_B \gamma_j \nabla u_j \cdot \nabla u_k \, \dd x - \int_{B \setminus \Omega} \gamma_j \nabla u_j \cdot \nabla u_k \, \dd x \\
= \int_B \gamma_j \nabla u_j \cdot \nabla (\gamma_j^{-1/2} v_k) \, \dd x + \int_B \gamma_j \nabla u_j \cdot \nabla ((\gamma_k^{-1/2} - \gamma_j^{-1/2}) v_k) \, \dd x\\
- \int_{B \setminus \Omega} \gamma_j \nabla u_j \cdot \nabla u_k \, \dd x.
\end{gather*}
Since $ u_j $ is a weak solution of $ \nabla \cdot (\gamma_j \nabla u_j) = 0 $ in $ \R^n $ and $ (\gamma_k^{-1/2} - \gamma_j^{-1/2}) v_k \in H^1_0(B) $, one has
\[ \int_B \gamma_j \nabla u_j \cdot \nabla ((\gamma_k^{-1/2} - \gamma_j^{-1/2}) v_k) \, \dd x = 0. \]
Thus,
\begin{equation*}
\duality{\Lambda_{\gamma_j|_\Omega}(u_j|_{\partial\Omega})}{u_k|_{\partial\Omega}} = \int_B \gamma_j \nabla (\gamma_j^{-1/2} v_j) \cdot \nabla (\gamma_j^{-1/2} v_k) \, \dd x - \int_{B \setminus \Omega} \gamma_j \nabla u_j \cdot \nabla u_k \, \dd x.
\end{equation*}
Using now that
\[ \duality{\Lambda_{\gamma_j|_\Omega} f}{g} = \duality{\Lambda_{\gamma_j|_\Omega} g}{f}, \]
for $ j \in \{ 1, 2 \} $, one gets
\begin{gather*}
\duality{(\Lambda_{\gamma_1|_\Omega} - \Lambda_{\gamma_2|_\Omegaº})(u_1|_{\partial\Omega})}{u_2|_{\partial\Omega}} + \int_{B \setminus \Omega} (\gamma_1 - \gamma_2) \nabla u_1 \cdot \nabla u_2 \, \dd x = \\
= \int_B \gamma_1 \nabla (\gamma_1^{-1/2} v_1) \cdot \nabla (\gamma_1^{-1/2} v_2) \, \dd x - \int_B \gamma_2 \nabla (\gamma_2^{-1/2} v_2) \cdot \nabla (\gamma_2^{-1/2} v_1) \, \dd x.
\end{gather*}
A simple computation shows
\begin{gather*}
\int_B \gamma_j \nabla (\gamma_j^{-1/2} v_j) \cdot \nabla (\gamma_j^{-1/2} v_k) \, \dd x = \int_B \nabla v_j \cdot \nabla v_k \, \dd x - \int_B \nabla \gamma_j^{1/2} \cdot \nabla (\gamma_j^{-1/2} v_j v_k) \, \dd x\\
= \int_B \nabla v_j \cdot \nabla v_k \, \dd x - \int_{\R^n} \nabla \gamma_j^{1/2} \cdot \nabla (\gamma_j^{-1/2} v_j v_k) \, \dd x.
\end{gather*}
Thus, we obtain the following Alessandrini formula
\begin{gather*}
\duality{(\Lambda_{\gamma_1|_\Omega} - \Lambda_{\gamma_2|_\Omega})(u_1|_{\partial\Omega})}{u_2|_{\partial\Omega}} + \int_{B \setminus \Omega} (\gamma_1 - \gamma_2) \nabla u_1 \cdot \nabla u_2 \, \dd x = \\
= \int_{\R^n} \nabla \gamma_2^{1/2} \cdot \nabla (\gamma_2^{-1/2} v_1 v_2) \, \dd x - \int_{\R^n} \nabla \gamma_1^{1/2} \cdot \nabla (\gamma_1^{-1/2} v_1 v_2) \, \dd x, \nonumber
\end{gather*}
which in turn implies \eqref{es:Integralbound}.
\end{proof}

This estimate has been proved for functions defined in $ \R^n $. Since this is not the case in the context of Calder\'on problem, we need to perform extensions of the coefficients. However, any kind of extension does not suffice since we need to control the term $ \norm{\gamma_1 - \gamma_2}{}{L^\infty(B \setminus \Omega)} $ on the right hand side of \eqref{es:Integralbound}. Thus, we are going to perform appropriate extensions of conductivities in $ C^{1,\varepsilon} (\overline{\Omega}) $ from their values on $ \partial \Omega $. These extensions are of Whitney type and they are based on the existence of certain polynomials approximating functions in $ C^{1,\varepsilon} (\overline{\Omega}) $. For the sake of completeness, we will show the existence of such polynomials in the next lemma. This makes necessary a quick explanation about what we mean by the space $ C^{1,\varepsilon} (\overline{\Omega}) $ and by $ \partial \Omega $ being locally described by the graph of a Lipschitz function.

We say that $ a \in C^{1,\varepsilon} (\overline{\Omega}) $, with $ 0 < \varepsilon \leq 1 $, if $ a : \Omega \longrightarrow \C $ is continuously differentiable in $ \Omega $ and its partial derivatives $ \partial^\alpha a $, with $ \alpha \in \N^n $ and $ |\alpha| \leq 1 $ satisfy
\begin{align}
& |\partial^\alpha a (x)| \leq C, \qquad \forall x\in \Omega, \quad |\alpha|\leq 1, \label{es:C1epsilonBOUNDNESS} \\
& |\partial^\alpha a (x) - \partial^\alpha a (y)| \leq C |x - y|^\varepsilon, \qquad \forall x,y \in \Omega, \quad |\alpha | = 1, \label{es:C1epsilonHOLDER}
\end{align}
for certain positive constant $ C $. The norm on $ C^{1,\varepsilon} (\overline{\Omega}) $, defined as the smallest constant $ C $ for which \eqref{es:C1epsilonBOUNDNESS} and \eqref{es:C1epsilonHOLDER} hold, makes $ C^{1,\varepsilon} (\overline{\Omega}) $ be a Banach space.

If $ \Omega $ is a bounded non-empty open subset of $ \R^n $, we say that $ \partial \Omega $ is locally described by the graph of a Lipschitz function if there exist $ \rho > 0 $, $ U_1, \dots, U_N $ open subsets of $ \R^n $, $ \y^1, \dots, \y^N $ isometric linear transformations in $ \R^n $ and $ \phi_1, \dots, \phi_N $ Lipschitz real-valued functions in $ \R^{n - 1} $ so that
\begin{itemize}
\item[(i)] if $ x^0 \in \partial \Omega $, then $ B(x^0; \rho) := \{ x \in \R^n : |x - x^0| < \rho \} \subset U_j $ for some $ j \in \{ 1, \dots, N \} $;
\item[(ii)] and $ \Omega \cap U_j = \{ x \in U_j : \phi_j (\y_1^j(x), \dots \y_{n - 1}^j(x)) > \y_n^j(x) \} $ for any $ j \in \{ 1, \dots, N \} $.
\end{itemize}

\begin{remark} \label{rem:parallelepiped} Since $ \Omega $ is bounded and $ \partial \Omega $ is locally described by the graph of a Lipschitz function, there exists a set $ \{ v^1, \dots, v^n \} $ of $ n $-linearly independent vectors in $ \R^n $ with $ \sum_{k=1}^n v^k = (0, \dots, 0, \rho') $ and $ \rho' > 0 $ such that, for any $ j \in \{ 1, \dots, N \} $,
\begin{equation*}
P (y) := \{ y + \sum_{k = 1}^n \lambda_k v^k : 0 < \lambda_k < 1, \} \subset \y^j (\Omega) \label{cond:inclusion}
\end{equation*}
for all $ y \in \y^j ( \Omega \cap V_j ) $ where $ V_j = \{ x \in U_j : \dist(x , U^\rho_j) < 2 \rho / 3 \} $ with $ U^\rho_j = \{ x \in U_j \cap \partial \Omega : B(x; \rho) \subset U_j \} $.

Furthermore, there exist $ r > 0 $ and $ c > 0 $ such that, for any $ j \in \{ 1, \dots, N \} $,
\begin{gather*}
Q_j(y^1, y^2) := P(y^1) \cap P(y^2) \neq \emptyset, \label{cond:nonEMPTY} \\
\inf_{z \in Q_j (y^1, y^2)} |y^1 - z| + |z - y^2| \leq c |y^1 - y^2|, \label{cond:equivalence}
\end{gather*}
for all $ y^1, y^2 \in \y^j ( \Omega \cap V_j ) $ such that $ |y^1 - y^2| < r $.
\end{remark}

Remark \ref{rem:parallelepiped} implies that, given any $ x^1 $ and $ x^2 $ in $ \Omega_j $ with $ \Omega_j $ being one of the connected components of $ \Omega $,  there exists a positive constant $ c $ (possibly different to $ c $ in Remark \ref{rem:parallelepiped}) such that
\[ \mathrm{dist}_{\Omega_j} (x^1, x^2) \leq c |x^1 - x^2|. \]
Here $ \mathrm{dist}_{\Omega_j} (x^1, x^2) $ is the infimum of the amounts $ \sum_{k = 0}^M |y^k - y^{k + 1}| $, where $ y^0 = x^1 $, $ y^{M + 1} = x^2 $ and $ t y^k + (1 - t) y^{k + 1} \in \Omega_j $ for all $ t \in [0, 1] $.

\begin{lemma} \sl If $ a \in C^{1, \varepsilon} (\overline{\Omega}) $, then $\partial^\alpha a $ with $ |\alpha| \leq 1 $ has a unique continuous extension $ f^{(\alpha)} : \overline{\Omega} \longrightarrow \C $ so that $\partial^{\alpha} a (x) = f^{(\alpha)} (x) $ for any $ x \in \Omega $ and if $ R_\alpha : \overline{\Omega} \times \overline{\Omega} \longrightarrow \C $ satisfies
\[ f^{(\alpha)}(x^1) = \sum_{|\beta + \alpha| \leq 1} f^{(\beta + \alpha)} (x^2) (x^1 - x^2)^\beta + R_\alpha (x^1, x^2), \qquad\forall x^1,x^2\in \overline{\Omega} \]
then
\begin{equation}
|R_\alpha (x^1, x^2)| \lesssim |x^1 - x^2|^{1 + \varepsilon - |\alpha|},\qquad\forall x^1,x^2\in \overline{\Omega}. \label{es:REMAINDERalpha}
\end{equation}
\end{lemma}

With these extensions of $ \partial^\alpha a $ to $ \overline{\Omega} $ in mind, we will make an abuse of notation identifying $ \partial^\alpha a $ with $ f^{(\alpha)} $ as functions defined in $ \overline{\Omega} $.

\begin{proof} Since $ \partial^\alpha a $ with $ |\alpha| = 1 $ is uniformly continuous and bounded in $ \Omega $, it can be uniquely extended to $ \overline{\Omega} $. Let $ f^{(\alpha)} $ denote its extension. 
On the other hand, since $ a $ is bounded and continuously differentiable in $ \Omega $ (whose boundary is locally described by the graph of a Lipschitz function), we know that $ a $ is a Lipschitz function in $ \Omega $. In particular, $ a $ is uniformly continuous and it can be extended to $ \overline{\Omega} $. Let $ f^{(0)} $ denote its extension.

The next thing to prove is \eqref{es:REMAINDERalpha}. The case $ |\alpha| = 1 $ is immediate. 
Note that in order to prove \eqref{es:REMAINDERalpha} for $ \alpha = 0 $, is enough to prove it for $ x^1 $ and $ x^2 $ in $ \Omega $. Reaching the whole $ \overline{\Omega} $ is a simple extension to the closure. Since $ \Omega $ is abounded Lipschitz domain there is a positive distance between their (finite number of) connected components. Thus, we only have to show that \eqref{es:REMAINDERalpha} holds for $ \alpha = 0 $ and $ x^1, x^2 \in \Omega_j $ with $ \Omega_j $ being any of the connected components of $ \Omega $. Let $ x^1 $ and $ x^2 $ be in $ \Omega_j $ and let $ y^1, \dots, y^M $ belong to $ \Omega_j $ such that $ t y^k + (1 - t) y^{k + 1} \in \Omega_j $ for all $ t \in [0, 1] $ and $ k \in \{ 0, \dots, M \} $ with $ y^0 = x^1 $ and $ y^{M + 1} = x^2 $. We have
\[ a(x^1) - a(x^2) = \sum_{|\alpha| = 1} \sum_{k = 0}^M (y^k - y^{k + 1})^\alpha \int_0^1 \partial^\alpha a (t y^k + (1 - t) y^{k + 1}) \, \dd t \]
and
\[ \sum_{|\alpha| = 1} \partial^\alpha a(x^2) (x^1 - x^2)^\alpha = \sum_{|\alpha| = 1} \sum_{k = 0}^M (y^k - y^{k + 1})^\alpha \int_0^1 \partial^\alpha a (y^{M + 1}) \, \dd t. \]
Therefore,
\[ |R_0(x^1, x^2)| \lesssim \sum_{k = 0}^M |y^k - y^{k + 1}|^{1 + \varepsilon}, \]
which in turn implies the result.
\end{proof}

Before carrying out the extension, we say that $ a \in C^{1, \varepsilon} (\R^n) $ if and only if it is bounded and continuously differentiable in $ \R^n $ and its partial derivatives $ \partial^\alpha a \in C^{0,\varepsilon} (\R^n) $ for $ |\alpha| = 1 $. Again, $ \norm{\centerdot}{}{C^{1,\varepsilon} (\R^n)} $ denotes the norm on $ C^{1, \varepsilon} (\R^n) $ defined as the minimum between $ \sup_{x \in \R^n} |a(x)| $ and $ \norm{\partial^\alpha a}{}{C^{0,\varepsilon} (\R^n)} $ with $ |\alpha| = 1 $.

\begin{lemma} \label{lem:extensionC1} \sl Consider a real constant $ \gamma_0 $ in the interval $ (0, 1] $. Let $ \gamma_1 $ and $ \gamma_2 $ belong to $ C^{1,\varepsilon}(\overline{\Omega}) $ such that $ \gamma_j (x) \geq \gamma_0 $ for all $ x \in \overline{\Omega} $ and $ j \in \{ 1, 2 \} $. There exist $ R > 0 $ and $ \sigma_1 $ and $ \sigma_2 $ in $ C^{1,\varepsilon}(\R^n) $ such that $ \overline{\Omega} \subset B $, $ \gamma_j = \sigma_j|_{\overline{\Omega}} $, $ \supp (\sigma_j - 1) \subset \overline{B} $, $ \sigma_j (x) \geq \gamma_0 / 2 $ for all $ x \in \R^n $,
\begin{equation}
\norm{\sigma_j}{}{C^{1,\varepsilon}(\R^n)} \lesssim \max \left( 1, \norm{\gamma_j}{}{C^{1,\varepsilon}(\overline{\Omega})} \right) \label{es:boundednessC1EX}
\end{equation}
and
\begin{equation}
\norm{\sigma_1 - \sigma_2}{}{L^\infty (\R^n \setminus \Omega)} \leq \sum_{|\alpha| \leq 1}\norm{\partial^\alpha \gamma_1 - \partial^\alpha \gamma_2}{}{L^\infty (\partial \Omega)}. \label{es:controlOUTSIDE}
\end{equation}
\end{lemma}

The implicit constant in \eqref{es:boundednessC1EX} only depends on $ n $.

\begin{proof} Let $ \epsilon_0 $ be a constant in the interval $ (0, 1) $ to be chosen later and set $ \Omega_0 = \{ x \in \R^n : \dist(x, \Omega) < \epsilon_0 \} $. Let $ F $ be the closed subset of $ \R^n $ given by $ F = \overline{\Omega} \sqcup (\R^n \setminus \Omega_0) $, where $\sqcup$ denotes disjoint union. Define $ f^{(\alpha)}_j : F \longrightarrow \R $ as $ f^{(\alpha)}_j (x) = \partial^\alpha \gamma_j (x) $ for any $ x \in \overline{\Omega} $ and $ |\alpha| \leq 1 $, $ f^{(0)}_j (x) = 1 $ for any $ x \in \R^n \setminus \Omega_0 $ and $ f^{(\alpha)}_j (x) = 0 $ for any $ x \in \R^n \setminus \Omega_0 $ and $ |\alpha| = 1 $. It is immediate that
\[ \left| f^{(\alpha)}_j (x) \right| < \max \left( 1, \norm{\gamma_j}{}{C^{1,\varepsilon}(\overline{\Omega})} \right), \qquad \forall x \in F,\, |\alpha| \leq 1  \]
and, if
\[ R_\alpha^j (x, y) = f^{(\alpha)}_j (x) - \sum_{|\alpha + \beta| \leq 1} f^{(\alpha + \beta)}_j(y) (x - y)^\beta  \qquad \forall x, y \in F,\, |\alpha| \leq 1 , \]
then
\begin{equation*}
|R_\alpha^j (x, y)| \lesssim |x - y|^{1 +\varepsilon - |\alpha|}
\end{equation*}
holds for all $|\alpha| \leq 1$ and all $ x, y \in F $.

In order to perform the proper extension, we follow Chapter VI, \S1 and \S2 in \cite{S}.  Write $ G = \R^n \setminus F $. It was proven in \cite{S} that one can define a partition of unity $ \{ \varphi^\ast_l \}_{l \in \N} $ subordinate to a collection of cubes $ \{ Q_l,\, Q_l^\ast \}_{l \in \N} $ with
\[ G = \bigcup_{l \in \N} Q_l^\ast,\qquad Q_l \subset Q_l^\ast \subset G\, \forall l \in \N, \]
satisfying
\[ \varphi^\ast_l (x) \geq 0, \, \sum_{l \in \N} \varphi^\ast_l (x) = 1\, \forall x \in G,\qquad Q_l \subset \mathrm{supp}\, \varphi^\ast_l \subset Q_l^\ast\, \forall l \in \N. \]
For any $ l \in \N $, choose a point $ y^l \in F $ such that $ \mathrm{dist}(Q_l, F) = \mathrm{dist}(Q_l, y^l) $. Note that $ y^l $ lays in $ \partial \Omega $ or in $ \partial \Omega_0 $. Now we are ready to present the extension given in \cite{S} (Chapter VI \S2). Define $ \sigma_j (x) = f^{(0)}_j (x) $ for any $ x \in F $ and
\[ \sigma_j (x) = \sum_{l \in \N} \left( \sum_{|\alpha| \leq 1} f^{(\alpha)}_j(y^l) (x - y^l)^\alpha \right) \varphi^\ast_l(x),\quad \forall x \in G. \]
In \cite{S} (Chapter VI \S2), it was stated that $ \sigma_j \in C^{1, \varepsilon}(\R^n) $ and \eqref{es:boundednessC1EX} holds. Furthermore, choosing $ R $ such that $ \Omega_0 \subset B $ one immediately has $ \supp (\sigma_j - 1) \subset \overline{B} $.

Let us next show that we can choose $ \epsilon_0 $ small enough such that $ \sigma_j (x) \geq \gamma_0 / 2 $. This last inequality immediately holds for $ x \in F $, so we just need to verify $ \sigma_j (x) \geq \gamma_0 / 2 $ for $ x \in G $. Note that for any $ x \in G $
\[ \sigma_j (x) \geq \gamma_0 - \sum_{|\alpha| = 1} \sup_{x \in \partial \Omega} |\partial^\alpha \gamma_j (x)| \sum_{l \in \N} |x - y^l| \varphi^\ast_l(x). \]
Before going any further, we recall some properties for the collections of cubes $ \{ Q_l,\, Q_l^\ast \}_{l \in \N} $ that were proven in Chapter VI \S1 of \cite{S}. These properties are the following: $ \mathrm{diam}(Q^\ast_l) \leq 5/4 \mathrm{diam}(Q_l) $ and $ \mathrm{diam}(Q_l) \leq \dist(F, Q_l) $. They imply that, for any $ x \in Q^\ast_l $
\[ |x - y^l| \leq \mathrm{diam}(Q^\ast_l) + \dist (Q_l, y^l) \leq 9/4 \dist(F, Q_l). \]
Since $ \dist(F, Q_l) \leq \epsilon_0 $ one gets
\[ \sigma_j (x) \geq \gamma_0 - \epsilon_0 \frac{9}{4} \sum_{|\alpha| = 1} \sup_{x \in \partial \Omega} |\partial^\alpha \gamma_j (x)|. \]
So, in order to get $ \sigma_j (x) \geq \gamma_0 / 2 $ for any $ x \in G $, it is enough to take $ \epsilon_0 \leq (2 / 9) \gamma_0 ( \sum_{|\alpha| = 1} \norm{\partial^\alpha \gamma_j}{}{L^\infty(\partial \Omega)} )^{-1} $.

Finally, \eqref{es:controlOUTSIDE} follows from
\[ |\sigma_1 (x) - \sigma_2 (x)| \leq \norm{\gamma_1 - \gamma_2}{}{L^\infty(\partial \Omega)} + \epsilon_0 \frac{9}{4} \sum_{|\alpha| = 1} \norm{\partial^\alpha \gamma_1 - \partial^\alpha \gamma_2}{}{L^\infty(\partial \Omega)} \qquad \forall x \in G \]
and the following choice
\[ \epsilon_0 \leq \min \left( \frac{4}{9}, \frac{2 \gamma_0}{9 \sum_{|\alpha| = 1} \norm{\partial^\alpha \gamma_1}{}{L^\infty(\partial \Omega)}}, \frac{2 \gamma_0}{9 \sum_{|\alpha| = 1} \norm{\partial^\alpha \gamma_2}{}{L^\infty(\partial \Omega)}} \right). \]
\end{proof}

We conclude this section formulating the estimate claimed at the beginning. It will be stated as a proposition.

\begin{proposition} \label{prop:extensionC1} \sl
Consider  a real constant $ \gamma_0 $ in the interval $ (0, 1] $. Let $ \gamma_1 $ and $ \gamma_2 $ belong to $ C^{1,\varepsilon} (\overline{\Omega}) $ such that $ \gamma_j (x) \geq \gamma_0 $ for all $ x \in \overline{\Omega} $ and $ j \in \{ 1, 2 \} $. There exist $ R > 0 $ and $ \sigma_1 $ and $ \sigma_2 $ in $ C^{1,\varepsilon} (\R^n) $ such that $ \overline{\Omega} \subset B $, $ \gamma_j = \sigma_j|_{\overline{\Omega}} $, $ \supp (\sigma_j - 1) \subset \overline{B} $, $ \sigma_j (x) \geq \gamma_0 / 2 $ for all $ x \in \R^n $,
\begin{equation*}
\sum_{|\alpha| \leq 1} \norm{\partial^\alpha \sigma_j}{}{L^\infty(\R^n)} \lesssim \max \left( 1, \sum_{|\alpha| \leq 1} \norm{\partial^\alpha \gamma_j}{}{L^\infty(\overline{\Omega})} \right) \qquad (\alpha \in \N^n)
\end{equation*}
and, for $ u_j \in H^1_\mathrm{loc}(\R^n) $ a weak solution of $ \nabla \cdot (\sigma_j \nabla u_j) = 0 $ in $ \R^n $, one has
\begin{gather}
\left| \int_{\R^n} \nabla \sigma_2^{1/2} \cdot \nabla (\sigma_2^{-1/2} v_1 v_2) \, \dd x - \int_{\R^n} \nabla \sigma_1^{1/2} \cdot \nabla (\sigma_1^{-1/2} v_1 v_2) \, \dd x \right| \lesssim \nonumber \\
\lesssim \left( \norm{\Lambda_{\gamma_1} - \Lambda_{\gamma_2}}{}{} + \norm{\Lambda_{\gamma_1} - \Lambda_{\gamma_2}}{\varepsilon/(1+\varepsilon)}{} \right) \norm{u_1}{}{H^1(B)} \norm{u_2}{}{H^1(B)}. \label{es:IntegralboundC1}
\end{gather}
Here $ v_j \in H^1_\mathrm{loc}(\R^n) $ denotes $ v_j = \sigma_j^{1/2} u_j $.
\end{proposition}

The implicit constant in \eqref{es:IntegralboundC1} depends on $ n, \Omega, \varepsilon, \gamma_0 $ and $ \norm{\gamma_j}{}{C^{1, \varepsilon} (\overline{\Omega})} $ for $ j = 1, 2 $.

\begin{proof}
From Lemma \ref{lem:integralBOUND} and \ref{lem:extensionC1} we immediately get
\begin{gather}
\left| \int_{\R^n} \nabla \sigma_2^{1/2} \cdot \nabla (\sigma_2^{-1/2} v_1 v_2) \, \dd x - \int_{\R^n} \nabla \sigma_1^{1/2} \cdot \nabla (\sigma_1^{-1/2} v_1 v_2) \, \dd x \right| \leq \nonumber \\
\leq \left( \norm{\Lambda_{\gamma_1} - \Lambda_{\gamma_2}}{}{} + \sum_{|\alpha| \leq 1}\norm{\partial^\alpha \gamma_1 - \partial^\alpha \gamma_2}{}{L^\infty (\partial \Omega)} \right) \norm{u_1}{}{H^1(B)} \norm{u_2}{}{H^1(B)}.
\end{gather}
On the other hand, it was proven by Alessandrini (see p. 256 in \cite{A90}) that
\[ \norm{\gamma_1 - \gamma_2}{}{L^\infty (\partial \Omega)} \lesssim \norm{\Lambda_{\gamma_1} - \Lambda_{\gamma_2}}{}{} \]
and
\[ \sum_{|\alpha| = 1}\norm{\partial^\alpha \gamma_1 - \partial^\alpha \gamma_2}{}{L^\infty (\partial \Omega)} \lesssim \norm{\Lambda_{\gamma_1} - \Lambda_{\gamma_2}}{\varepsilon/(1+\varepsilon)}{}, \]
where the implicit constants depend on $ n, \Omega, \varepsilon, \gamma_0 $ and $ \norm{\gamma_j}{}{C^{1, \varepsilon} (\overline{\Omega})} $ for $ j = 1, 2 $. This proves the statement of this proposition.
\end{proof}

\section{Complex geometrical optics} \label{sec:CGOreview}
In this section we review the construction of CGOs for the conductivity equation following the arguments presented in \cite{HT}. Here, CGOs are weak solutions of $ \nabla \cdot (\gamma \nabla u) = 0 $ in $ \R^n $ having the special form:
\[ u_\zeta = \gamma^{-1/2} e^\rho (1 + r_\zeta) \]
with $ \rho(x) = \zeta \cdot x $ and $ \zeta \in \C^n $ satisfying that $ \zeta \cdot \zeta = 0 $ and $ |\zeta| \geq 1 $. Along these notes $ \zeta_1 \cdot \zeta_2 $ with $ \zeta_j \in \C^n $ stands for the analytic extension of the real-inner product. Here $ r_\zeta $ has to be understood as a remainder going to zero in some sense. Following the ideas in \cite{SU} and \cite{B}, one can transform the conductivity equation into a Schr\"odinger equation by \textit{rescaling} the solution. More precisely, if $ \gamma \in C^{0,1}(\R^n) $ one can check that $ u \in H^1_\mathrm{loc}(\R^n) $ is a weak solution of $ \nabla \cdot (\gamma \nabla u) = 0 $ in $ \R^n $ if and only if $ v = \gamma^{1/2} u \in H^1_\mathrm{loc}(\R^n) $ satisfies
\[ \int_{\R^n} \nabla v \cdot \nabla \varphi \, \dd x - \int_{\R^n} \nabla \gamma^{1/2} \cdot \nabla (\gamma^{-1/2} v \varphi) \, \dd x = 0 \]
for any $ \varphi \in C^\infty_0(\R^n) $ --the space of smooth functions with compact support. Following the notation used in \cite{SU} and \cite{B}, we shall write
\begin{align}
\duality{q}{\varphi} &:= - \int_{\R^n} \nabla \gamma^{1/2} \cdot \nabla (\gamma^{-1/2} \varphi) \, \dd x \label{def:q}, \\
\duality{m_q v}{\varphi} &:= - \int_{\R^n} \nabla \gamma^{1/2} \cdot \nabla (\gamma^{-1/2} v \varphi) \, \dd x, \label{def:mq}
\end{align}
for any $ \varphi \in C^\infty_0(\R^n) $. Thus, in order to construct CGOs for the conductivity equation it is enough to prove the existence of $ r_\zeta $ satisfying
\begin{equation}\label{eq:remainder}
- \Delta r_\zeta - 2 \zeta \cdot \nabla r_\zeta + m_q r_\zeta = - q
\end{equation}
in $ \R^n $ and to deduce the reminder properties for $ r_\zeta $. This task has been carried out in \cite{SU} for smooth conductivities, in \cite{B} for conductivities in $ \cup_{\varepsilon > 0} C^{1,1/2 + \varepsilon}(\overline{\Omega}) $, and in \cite{HT} for continuously differentiable conductivities and small enough Lipschitz conductivities.

Haberman and Tataru introduced in \cite{HT} a space similar to Bourgain's $ X^{s,b} $-spaces to study equation (\ref{eq:remainder}) for the remainder. Concretely, they introduced the homogeneous space $ \dot{X}_\zeta^b $ defined as follows, $ u \in \dot{X}_\zeta^b $ if and only if $ u \in \mathcal{S}'(\R^n) $ (the space of tempered distributions) and $ \widehat{u} $, the Fourier transform of $ u $, belongs to $ L^2(\R^n; |p_\zeta|^{2b} \dd \xi) $ (the $ L^2 $-space in $ \R^n $ with respect to the measure $ |p_\zeta|^{2b} \dd \xi $). Here $ b \in \R $ and $ p_\zeta (\xi) := |\xi|^2 - 2i \zeta \cdot \xi $ (the symbol of the conjugated laplacian $ - \Delta_\zeta := - \Delta - 2 \zeta \cdot \nabla $). This space endowed with the norm
\[ \norm{u}{}{\dot{X}_\zeta^b} := \norm{|p_\zeta|^b \widehat{u}}{}{L^2(\R^n)} \]
is Banach space at least when $ b < 1 $. Note that the operator norm of $ (- \Delta_\zeta)^{-1} : \dot{X}_\zeta^{-1/2} \longrightarrow \dot{X}_\zeta^{1/2} $ (defined by the symbol $ 1 / p_\zeta $) is
\[ \norm{(- \Delta_\zeta)^{-1}}{}{\mathcal{L}\left( \dot{X}_\zeta^{-1/2}, \dot{X}_\zeta^{1/2} \right)} = 1. \]
Note that whenever
\begin{equation}\label{es:mq}
\norm{m_q}{}{\mathcal{L}\left( \dot{X}_\zeta^{1/2}, \dot{X}_\zeta^{-1/2} \right)} < 1,
\end{equation}
one knows that
\[ r_\zeta = \left( I + (- \Delta_\zeta)^{-1} m_q \right)^{-1} \left( (- \Delta_\zeta)^{-1}(- q) \right) \]
is a solution of (\ref{eq:remainder}), just by the Neumann series. Additionally, we can estimate the remainder $ r_\zeta $ in \eqref{eq:remainder} by
\begin{equation}\label{ter:normq}
\norm{q}{}{\dot{X}_\zeta^{-1/2}}.
\end{equation}
Indeed, in that case we would have
\begin{equation}
\norm{r_\zeta}{}{\dot{X}_\zeta^{1/2}} \leq \norm{(I + (- \Delta_\zeta)^{-1} m_q)^{-1}}{}{\mathcal{L}\left( \dot{X}_\zeta^{1/2} \right)} \norm{q}{}{\dot{X}_\zeta^{-1/2}}. \label{es:rBOUNDBYq}
\end{equation}

Let us now sketch how Haberman and Tataru proved (\ref{es:mq}) in the case where $ \gamma $ is continuously differentiable in $ \R^n $ and constant outside $ B $. Afterwards, we show how they deduced the remainder properties of the $ r_\zeta $. By the Leibniz rule, we have that
\begin{gather}
\duality{m_q u}{v} 
=- \int_{\R^n} \nabla \gamma^{1/2} \cdot \nabla \gamma^{-1/2} u v \,\dd x - \int_{\R^n} \gamma^{-1/2} \nabla \gamma^{1/2} \cdot \nabla(uv) \,\dd x. \label{eq:Leibniz}
\end{gather}
The first term on the right hand side of \eqref{eq:Leibniz} was estimated by
\begin{equation}\label{es:1termLEIBNIZ}
\left|\int_{\R^n}\nabla\gamma^{1/2}\cdot\nabla\gamma^{-1/2}uv\,\dd x\right|\lesssim \sum_{j = 1}^n \norm{\partial_{x_j} \log \gamma}{2}{L^\infty(\R^n)} |\zeta|^{-1}\|u\|_{\BHP}\|v\|_{\BHP}
\end{equation}
in \cite{HT} (Corollary 2.1)\footnote{Along this section the implicit constants only depend on $ n $ and $ \Omega $.}
. In order to prove that the second term on the right hand side of \eqref{eq:Leibniz} allows \eqref{es:mq} to be held, Haberman and Tataru required $ \nabla \log \gamma $ to be of compact support as well as in $ C^0(\R^n) $ to eventually use an approximation to the identity. Here $ C^0(\R^n) $ denotes the space of continuous functions in $ \R^n $. Define $ \psi_{h}(x) = h^{-n} \psi (x / h) $, where $ h > 0 $ and $ \psi $ is a smooth function in $ \R^n $ supported on the unit ball satisfying $ \int_{\R^n} \psi \, \dd x = 1 $, and write
\begin{gather*}
\left| \int_{\R^n} \gamma^{-1/2} \nabla \gamma^{1/2} \cdot \nabla(uv) \, \dd x \right| \leq \left| \int_{\R^n} \psi_{h} \ast (\gamma^{-1/2} \nabla \gamma^{1/2}) \cdot \nabla(uv) \, \dd x \right| \\
+\left| \int_{\R^n} (\psi_h \ast (\gamma^{-1/2} \nabla \gamma^{1/2}) - \gamma^{-1/2} \nabla \gamma^{1/2}) \cdot \nabla(uv) \, \dd x \right|.
\end{gather*}
In Lemma 2.3 of \cite{HT} the following estimates were proven:
\begin{equation*}
\left| \int_{\R^n} \psi_{h} \ast (\gamma^{-1/2} \nabla \gamma^{1/2}) \cdot \nabla(uv) \, \dd x \right| \lesssim \frac{1}{h |\zeta|} \sum_{j = 1}^n \norm{\partial_{x_j} \log \gamma}{}{L^\infty(\R^n)} \norm{u}{}{\BHP} \norm{v}{}{\BHP}
\end{equation*}
and
\begin{gather*}
\left| \int_{\R^n}(\psi_h \ast (\gamma^{-1/2} \nabla\gamma^{1/2}) - \gamma^{-1/2} \nabla \gamma^{1/2}) \cdot \nabla(uv) \, \dd x \right| \lesssim \\
\lesssim \sum_{j = 1}^n \norm{\psi_h \ast (\gamma^{-1/2} \partial_{x_j} \gamma^{1/2}) - \gamma^{-1/2} \partial_{x_j} \gamma^{1/2})}{}{L^\infty(\R^n)} \norm{u}{}{\BHP} \norm{v}{}{\BHP}.
\end{gather*}
Before proceeding let us introduce some notation. Let $ f $ be either in $ L^p(\R^n) $ with $ 1 \leq p < + \infty $ or in $ L^\infty(\R^n) \cap C^0(\R^n) $ and define the $ L^p $-modulus of continuity as
\begin{equation}\label{for:modulusCONT}
\omega_p f(t) := \sup_{|y| < 1} \norm{f - f(\centerdot - ty)}{}{L^p(\R^n)}
\end{equation}
with $ p \in [1, +\infty) $ or $ p = \infty $. Note that
\[ \norm{\psi_h \ast (\gamma^{-1/2} \partial_{x_j} \gamma^{1/2}) - \gamma^{-1/2} \partial_{x_j} \gamma^{1/2})}{}{L^\infty(\R^n)} \leq \omega_\infty (\partial_{x_j} \log \gamma)(h). \]
We now take $ h = |\zeta|^{-1/(1 + \varepsilon)}$ and check that
\begin{gather}
\norm{m_q}{}{\mathcal{L}\left( \dot{X}_\zeta^{1/2}, \dot{X}_\zeta^{-1/2} \right)} \lesssim  |\zeta|^{-1} \sum_{j = 1}^n \norm{\partial_{x_j} \log \gamma}{2}{L^\infty(\R^n)} \label{es:mqCHOOSEepsilon} \\
+ |\zeta|^{-\varepsilon / (1 + \varepsilon)} \sum_{j = 1}^n \norm{\partial_{x_j} \log \gamma}{}{L^\infty(\R^n)} + \sum_{j = 1}^n \omega_\infty (\partial_{x_j}\log \gamma)(|\zeta|^{-1/(1 + \varepsilon)}). \nonumber
\end{gather}
Hence (\ref{es:mq}) holds for $ |\zeta| $ large enough.

The next step will be to deduce the remainder properties of the solution of (\ref{eq:remainder}). As we pointed out previously, it will be enough to study (\ref{ter:normq}) and to check if this norm tends to vanish in some sense. Whenever the conductivity is smooth enough (for instance being in $ H^{3/2 + \varepsilon}(\R^n) $) one can prove that
\[ \norm{q}{}{\dot{X}_\zeta^{-1/2}} \lesssim |\zeta|^{-\varepsilon}, \]
for any $ \varepsilon \in (0, 1/2] $. However, this kind of estimates seems to fail for less regular conductivities. Despite this, Haberman and Tataru showed that (\ref{ter:normq}) decays in average for some choices of $ \zeta $.

Firstly note that
\begin{equation}\label{eq:q}
\duality{q}{v} 
=\frac{1}{4} \int_{\R^n} |\nabla \log \gamma|^2 v \, \dd x - \frac{1}{2} \int_{\R^n} \nabla \log \gamma \cdot \nabla  v \, \dd x.
\end{equation}
In order to estimate the first term in \eqref{eq:q}, we note that
\begin{equation*}
\int_{\R^n} |\nabla \log \gamma|^2 v \, \dd x = \int_{\R^n} |\nabla \log \gamma|^2 \phi v \, \dd x,
\end{equation*}
for any radial function $ \phi \in C^\infty_0(\R^n) $ such that $ \phi(x) = 1 $ for all $ x \in B $. From now on, $ \phi $ will only denote a function with these properties. Hence
\[ \left| \int_{\R^n} |\nabla \log \gamma|^2 v \, \dd x \right| \lesssim \sum_{j = 1}^n \norm{\partial_{x_j} \log \gamma}{2}{L^\infty(\R^n)} \norm{\phi v}{}{L^2(\R^n)}. \]
Using Lemma 2.2 from \cite{HT}, one gets
\[ \left| \int_{\R^n} |\nabla \log \gamma|^2 v \, \dd x \right| \lesssim |\zeta|^{-1/2} \sum_{j = 1}^n \norm{\partial_{x_j} \log \gamma}{2}{L^\infty(\R^n)} \norm{v}{}{\dot{X}^{1/2}_\zeta}. \]
For the second term of \eqref{eq:q}, it also holds
\[ - \int_{\R^n} \nabla \log \gamma \cdot \nabla v \, \dd x = \duality{\phi\, \Delta \log \gamma}{v}. \]
Hence
\begin{equation*}
\left| \int_{\R^n} \nabla \log \gamma \cdot \nabla v \, \dd x \right| \leq \|\phi\, \Delta \log\gamma \|_{\BHM} \|v\|_{\BHP}.
\end{equation*}
Therefore,
\begin{equation}\label{eq:boundq}
\|q\|_{\BHM} \lesssim |\zeta|^{-1/2} \sum_{j = 1}^n \norm{\partial_{x_j} \log \gamma}{2}{L^\infty(\R^n)} + \|\phi\, \Delta \log\gamma \|_{\BHM}.
\end{equation}
Obviously, the second term in \eqref{eq:boundq} will determine the properties of \eqref{ter:normq}. Actually, if one just does some straight computations one only gets:
\begin{gather}
\norm{\phi\, \Delta \log\gamma}{}{\BHM} \leq \norm{\Delta \log\gamma}{}{\BM} \lesssim \nonumber \\
\lesssim \left( \int_{|\xi| < 4 |\zeta|} |\zeta| |\widehat{\nabla \log\gamma}(\xi)|^2 \, \dd \xi + \int_{|\xi| \geq 4 |\zeta|} |\widehat{\nabla \log\gamma}(\xi)|^2 \, \dd \xi \right)^{1/2} \nonumber \\
\lesssim (1 + |\zeta|)^{1/2} \left( \sum_{j = 1}^n \norm{\partial_{x_j} \log \gamma}{2}{L^\infty(\R^n)} \right)^{1/2}, \label{es:NONdecay}
\end{gather}
where we used Lemma 2.2 from \cite{HT}, with
\[ \norm{u}{}{X_\zeta^b} := \norm{(|\zeta| + |p_\zeta|)^b \widehat{u}}{}{L^2(\R^n)} \]
for $ b \in \R $, and $ |\xi|^2/2 \leq |p_\zeta(\xi)| \leq 3|\xi|^2/2 $ whenever $ 4 |\zeta| \leq |\xi| $. Since this term does not decay as $ |\zeta| $ grows for low regular conductivities
, Haberman and Tataru studied its behaviour \textit{in average}. They first took an arbitrary $ k \in \R^n $, consider $ P $ a $ 2 $-dimensional linear subspace orthogonal to $ k $ and set
\begin{equation}
\zeta := s \eta + i \left( \frac{k}{2} + \left( s^2 - \frac{|k|^2}{4} \right)^{1/2} \kappa \right), \label{def:vectorZETA}
\end{equation}
where $ s \in [|k|/2, +\infty) $, $ \eta \in P \cap \{ x \in \R^n : |x| = 1 \} $ (for later references set $ S := P \cap \{ x \in \R^n : |x| = 1 \} $) and $ \kappa $ is the unique vector making $ \{ \eta, \kappa \} $ a positively oriented orthonormal basis of $ P $. This kind of choice for $ \zeta $ became standard after \cite{C} and \cite{SU}. Note that $ \zeta \cdot \zeta = 0 $, $ |\zeta|^2 = 2 s^2 $ and $ S $ depends on $ |k|^{-1} k \in \{ x \in \R^n : |x| = 1 \} $. For a fix $ k \in \R^n $, $ \zeta $ only depends on $ s \in [|k|/2, +\infty) $ and $ \eta \in S $. At this point, Haberman and Tataru integrated the function
\[ (s, \eta) \in [\lambda, 2\lambda] \times S \longmapsto \|\phi\, \Delta \log \gamma \|^2_{\BHM} \]
with respect to $ \lambda^{-1} \dd s\, \dd l $, where $ \dd l $ stands for the length form on $ S $. Thus,
\begin{gather*}
\frac{1}{\lambda} \int_S \int_{\lambda}^{2\lambda} \|\phi\, \Delta \log \gamma \|^2_{\BHM}\, \dd s \, \dd l \lesssim \frac{1}{\lambda} \int_S \int_{\lambda}^{2\lambda} \|\phi\, \nabla \cdot (\psi_h \ast \nabla \log \gamma) \|^2_{\BHM}\, \dd s \, \dd l \\ + \frac{1}{\lambda} \int_S \int_{\lambda}^{2\lambda} \|\phi\, \nabla \cdot (\psi_h \ast \nabla \log \gamma - \nabla \log \gamma) \|^2_{\BHM}\, \dd s \, \dd l,
\end{gather*}
where $ \lambda \geq |k| / 2 $ and $ \psi_h $ stands for the approximation to the identity previously introduced. In Lemma 3.1 of \cite{HT}, each term on the right hand side was estimated in such a way that
\begin{gather*}
\frac{1}{\lambda} \int_S \int_{\lambda}^{2\lambda} \|\phi\, \Delta \log \gamma \|^2_{\BHM}\, \dd s \, \dd l \lesssim \frac{1}{\lambda}  \norm{\nabla \cdot (\psi_h \ast \nabla \log \gamma)}{2}{L^2(\R^n)} \\
+ (1+\langle k\rangle^2 / \lambda)\sum_{j = 1}^n \norm{\psi_h \ast \partial_{x_j} \log \gamma - \partial_{x_j} \log \gamma}{2}{L^2(\R^n)},
\end{gather*}
where $\langle k\rangle=(1+|k|^2)^{1/2}$ and $ |k| \leq \lambda $. On one hand,
\[ \norm{\nabla \cdot (\psi_h \ast \nabla \log \gamma)}{2}{L^2(\R^n)} \lesssim \frac{1}{h^2} \sum_{j = 1}^n\norm{\partial_{x_j} \log \gamma}{2}{L^2(\R^n)}. \]
On the other hand,
\[ \norm{\psi_h \ast \partial_{x_j} \log \gamma - \partial_{x_j} \log \gamma}{}{L^2(\R^n)} \leq \omega_2 (\partial_{x_j} \log \gamma)(h). \]
The notation $ \omega_p $ was introduced in \eqref{for:modulusCONT}. Therefore
\begin{gather}
\frac{1}{\lambda} \int_S \int_{\lambda}^{2\lambda} \|\phi\, \Delta \log \gamma \|^2_{\BHM}\, \dd s \, \dd l \lesssim \frac{1}{h^2 \lambda} \sum_{j = 1}^n \norm{\partial_{x_j} \log \gamma}{2}{L^2(\R^n)} \nonumber \\
+ (1+\langle k\rangle^2 / \lambda)\sum_{j = 1}^n \big( \omega_2 (\partial_{x_j} \log \gamma)(h) \big)^2, \label{es:averageDECAYlaplacian}
\end{gather}
with $ |k| \leq \lambda $. As we did in \eqref{es:mqCHOOSEepsilon}, we could now choose $ h $ as a negative power of $ \lambda $, however, in order to make the optimal choice for our case, we will wait until Section \ref{sec:stability}. Thus, we have by now that, for $ \lambda \geq |k| $,
\begin{gather}
\frac{1}{\lambda} \int_S \int_{\lambda}^{2\lambda} \|q\|^2_{\BHM}\, \dd s \, \dd l \lesssim \frac{1}{\lambda} \left( \sum_{j = 1}^n \norm{\partial_{x_j} \log \gamma}{2}{L^\infty(\R^n)} \right)^2 \label{eq:potentiallarge} \\
+ \frac{1}{h^2 \lambda} \sum_{j = 1}^n \norm{\partial_{x_j} \log \gamma}{2}{L^2(\R^n)} + (1+\langle k\rangle^2 / \lambda)\sum_{j = 1}^n \big( \omega_2 (\partial_{x_j} \log \gamma)(h) \big)^2. \nonumber
\end{gather}

We end this section by computing the $H^1(B)$-norm of the CGOs constructed by Haberman and Tataru and proving the estimate
\begin{equation}\label{es:multiBYphi}
\norm{e^{i k \cdot x} \phi \, w}{}{X^{1/2}_\zeta} \lesssim \langle k \rangle^{1/2} \norm{w}{}{\dot{X}^{1/2}_\zeta}
\end{equation}
for any $ w \in \dot{X}^{1/2}_\zeta $, for any $ k \in \R^n $ and $ \zeta $ as in \eqref{def:vectorZETA} with $ s \geq 1 $. The implicit constant in this estimate only depends on $ \phi $. This estimate was essentially proven in Lemma 2.2 in \cite{HT} but Haberman and Tataru did not make explicit there the dependence on $ k $ of the implicit constant.\footnote{In the first version of our manuscript, we stated that the constant of estimate \eqref{es:multiBYphi} did not depend on $ k $. This was a mistake that was pointed out by the anonymous referee. We also thank him or her for this.}

Let us first compute the $H^1(B)$-norm of the CGOs. By definition
\begin{equation}\label{eq:H1norm}
\|u_{\zeta}\|^2_{H^1(B)}=\|u_{\zeta}\|^2_{L^2(B)}+\sum_{j = 1}^n\|\partial_{x_j}u_{\zeta}\|^2_{L^2(B)}.
\end{equation}
Concerning the first term in \eqref{eq:H1norm} we get
\begin{equation*}
\|u_{\zeta}\|_{L^2(B)} \leq \gamma_0^{-1/2} e^{R|\zeta|} \left( R^{n/2} + \|r_{\zeta}\|_{L^2(B)} \right).
\end{equation*}
Furthermore, using the same $ \phi $ previously introduced one sees
\begin{equation}\label{es:L2toX1/2}
\|r_{\zeta}\|_{L^2(B)} \leq \norm{\phi r_\zeta}{}{L^2(\R^n)} \leq |\zeta|^{-1/2} \|\phi r_\zeta\|_{\BP}\lesssim |\zeta|^{-1/2}\|r_{\zeta}\|_{\BHP},
\end{equation}
where the last inequality follows by Lemma 2.2 in \cite{HT}. On the other hand, the first derivative terms in \eqref{eq:H1norm} can be bounded as follows
\begin{gather*}
\|\partial_{x_j}u_{\zeta}\|_{L^2(B)}  \leq \gamma_0^{-1/2} e^{R|\zeta|} \left( \|\partial_{x_j}\log\gamma\|_{L^{\infty}(\R^n)} + |\zeta_j| \right) \left( R^{n/2} + \|r_{\zeta}\|_{L^2(B)} \right) \\
+ \gamma_0^{-1/2} e^{R|\zeta|} \norm{\partial_{x_j} r_\zeta}{}{L^2(B)}.
\end{gather*}
We have again
\begin{gather*}
\sum_{j = 1}^n \|\partial_{x_j} r_{\zeta}\|^2_{L^2(B)} \leq \sum_{j = 1}^n \norm{\partial_{x_j}(\phi r_\zeta)}{2}{L^2(\R^n)} = \int_{\R^n} |\xi|^2 |\widehat{\phi r_\zeta}(\xi)|^2 \, \dd \xi \\
\lesssim \int_{|\xi| < 4 |\zeta|} |\zeta|^2 |\widehat{\phi r_\zeta}(\xi)|^2 \, \dd \xi + \int_{|\xi| \geq 4 |\zeta|} |p_\zeta(\xi)| |\widehat{\phi r_\zeta}(\xi)|^2 \, \dd \xi \\
\lesssim |\zeta|^2 \norm{\phi r_\zeta}{2}{L^2(\R^n)} + \norm{\phi r_\zeta}{2}{X^{1/2}_\zeta} \\
\lesssim (|\zeta| + 1) \norm{r_\zeta}{2}{\dot{X}^{1/2}_\zeta}.
\end{gather*}
Here we used again that $ |\xi|^2/2 \leq |p_\zeta(\xi)| \leq 3|\xi|^2/2 $ whenever $ 4 |\zeta| \leq |\xi| $, estimate (\ref{es:L2toX1/2}) and Lemma 2.2 from \cite{HT}. Summing up,
\begin{gather}
\norm{u_\zeta}{}{H^1(B)} \lesssim \gamma_0^{-1/2} e^{R|\zeta|} \left( \sum_{j = 1}^n \|\partial_{x_j}\log\gamma\|_{L^{\infty}(\R^n)} + |\zeta| \right) \left( R^{n/2}+ |\zeta|^{-1/2}\|r_{\zeta}\|_{\BHP} \right) \nonumber \\
+ \gamma_0^{-1/2} e^{R|\zeta|} (1 + |\zeta|)^{1/2} \norm{r_\zeta}{}{\dot{X}^{1/2}_\zeta}. \label{eq:normaH1}
\end{gather}

Finally, let us prove estimate \eqref{es:multiBYphi}. Let $ \phi_k $ denote the function $ \phi_k(x) = e^{i k \cdot x} \phi(x) $ for any $ x \in \R^n $. Since
\[ \widehat{\phi_k w} = (2 \pi)^{-n} \widehat{\phi_k} \ast \widehat{w}, \]
we see that, in order to prove estimate \eqref{es:multiBYphi}, it is enough to show that the convolution operator
\[ v \in L^2(\R^n; |p_\zeta| d\xi) \longmapsto \widehat{\phi_k} \ast v \in L^2(\R^n; (|\zeta| + |p_\zeta|) \, d\xi) \]
is bounded and its norm is bounded by a multiple of $ \langle k \rangle^{1/2} $. In turn, by Lemma 2.1 in \cite{HT} we only have to check that
\[ \sup_{x \in \R^n} \int_{\R^n} \frac{|\zeta| + |p_\zeta(x)|}{|p_\zeta(y)|} |\widehat{\phi_{k}} (x - y)| \, \dd y \lesssim \langle k \rangle. \]
Note that
\[ \int_{\R^n} \frac{|\zeta| + |p_\zeta(x)|}{|p_\zeta(y)|} |\widehat{\phi_{k}} (x-y)| \, \dd y = \int_{\R^n} \frac{|\zeta| + |p_\zeta(x)|}{|p_\zeta(y - k)|} |\widehat{\phi} (y - x)| \, \dd y \]
for any $ x \in \R^n $ and $ |p_\zeta (y - k)| = |p_{\zeta'} (y)| $ for any $ y, k \in \R^n $, where
\[ \zeta' := s \eta + i \left(- \frac{k}{2} + \left( s^2 - \frac{|k|^2}{4} \right)^{1/2} \kappa \right). \]
Since $ |\zeta|^2 = 2 s^2 $, we will study the integral
\begin{gather}
\int_{\R^n} \frac{s + |p_\zeta(x)|}{|p_{\zeta'}(y)|} |\widehat{\phi} (y - x)| \, \dd y \lesssim \nonumber \\
\lesssim \int_{|y| \geq 8s} \frac{s + |p_\zeta(x)|}{|y|^2} |\widehat{\phi} (y - x)| \, \dd y + \int_{|y| < 8s} \frac{s + |p_\zeta(x)|}{|p_{\zeta'}(y)|} |\widehat{\phi} (y - x)| \, \dd y. \label{es:lowShightS}
\end{gather}
Here we used that $ |y|^2 \leq 2 |p_{\zeta'}(y)| $ when $ |y| \geq 8s $. Let us start by the first integral in the right hand side of \eqref{es:lowShightS}. Assume $ |x| \geq 4 |\zeta| $, then $ |p_\zeta (x)| \leq 3 ( |x - y|^2 + |y|^2 ) $ and
\[ \int_{|y| \geq 8s} \frac{s + |p_\zeta(x)|}{|y|^2} |\widehat{\phi} (y - x)| \, \dd y \lesssim \int_{|y| \geq 8s} \frac{s + |x - y|^2 + |y|^2}{|y|^2} |\widehat{\phi} (y - x)| \, \dd y \lesssim 1. \]
Assume now $ |x| < 4 |\zeta| $, then $ |p_\zeta (x)| \leq 24 |\zeta|^2 $ and
\[ \int_{|y| \geq 8s} \frac{s + |p_\zeta(x)|}{|y|^2} |\widehat{\phi} (y - x)| \, \dd y \lesssim \int_{|y| \geq 8s} \frac{s + s^2}{s^2} |\widehat{\phi} (y - x)| \, \dd y \lesssim 1. \]
Let us continue by the second integral in the right hand side of \eqref{es:lowShightS}. Note that $ |y| < 8s $ implies
\[ |p_{\zeta'} (y)| \sim s ( |s - |y + \mathrm{Im}\,\zeta'|| + |y \cdot \eta| ). \]
Assume $ |x| \geq 8 |\zeta| $, then $ |p_\zeta (x)| \leq 3 ( |x - y|^2 + |y|^2 ) $, $ (\sqrt{2} - 1) 8 s < |x| - |y| \leq |x - y| $ and
\[ \int_{|y| < 8s} \frac{s + |p_\zeta(x)|}{|p_{\zeta'}(y)|} |\widehat{\phi} (y - x)| \, \dd y \lesssim \int_{|y| < 8s} \frac{(1 + |x - y|^2)|\widehat{\phi} (y - x)|}{|s - |y + \mathrm{Im}\,\zeta'|| + |y \cdot \eta|}  \, \dd y \lesssim 1. \]
Last estimate is a consequence of (8) in the proof of Lemma 2.2 of \cite{HT}. Finally, assume $ |x| < 8 |\zeta| $, then
\[ |p_{\zeta} (x)| \sim s ( |s - |x + \mathrm{Im}\,\zeta|| + |x \cdot \eta| ) \]
and
\[ \int_{|y| < 8s} \frac{s + |p_\zeta(x)|}{|p_{\zeta'}(y)|} |\widehat{\phi} (y - x)| \, \dd y \lesssim \int_{|y| < 8s} \frac{1 + |s - |x + \mathrm{Im}\,\zeta|| + |x \cdot \eta|}{|s - |y + \mathrm{Im}\,\zeta'|| + |y \cdot \eta|} |\widehat{\phi} (y - x)| \, \dd y. \]
By the definitions of $ \zeta $ and $ \zeta' $, we have
\begin{equation}
|s - |x + \mathrm{Im}\,\zeta|| \leq |x - y| + |k| + |s - |y + \mathrm{Im}\,\zeta'||, \qquad |x \cdot \eta| \leq |x - y| + |y \cdot \eta|. \label{es:diffxy}
\end{equation}
Using \eqref{es:diffxy} and (8) in the proof of Lemma 2.2 of \cite{HT}, we have
\[ \int_{|y| < 8s} \frac{s + |p_\zeta(x)|}{|p_{\zeta'}(y)|} |\widehat{\phi} (y - x)| \, \dd y \lesssim 1 + |k|. \]
This concludes the proof of \eqref{es:multiBYphi}.

\section{Stability estimates}\label{sec:stability}
Along this section, $ \gamma_1 $ and $ \gamma_2 $ will be in $ C^{1,\varepsilon}(\overline{\Omega}) $ such that $ \gamma_j (x) > 1/M $ for all $ x \in \overline{\Omega} $ and $ j \in \{ 1, 2 \} $. Furthermore, we shall assume
\[ \norm{\gamma_j}{}{C^{1,\varepsilon}(\overline{\Omega})} \leq M. \]
We first extend $ \gamma_1 $ and $ \gamma_2 $ as in Lemma \ref{lem:extensionC1} (note that $ \epsilon_0 $ defining $ \Omega_0 $ in the proof of Lemma \ref{lem:extensionC1} is the same for any $ \gamma_1 $ and $ \gamma_2 $ in these conditions). Let their extensions be denoted by $ \gamma_j $ instead of $ \sigma_j $. Consider $ q_j $ with $ j \in \{ 1, 2 \} $ defined as in \eqref{def:q} for $ \gamma_j $ instead of $ \gamma $. That definition can be extended to $ \varphi \in H^1(\R^n) $ which implies that $ q_j \in H^{-1}(\R^n) $. Moreover,
\begin{gather*}
\duality{q_1 - q_2}{\varphi} = - \int_{\R^n} \nabla( \log \gamma_1^{1/2} - \log \gamma_2^{1/2} ) \cdot \nabla \varphi \, \dd x \\
+ \int_{\R^n} \nabla \log( \gamma_1^{1/2} \gamma_2^{1/2} ) \cdot \nabla ( \log \gamma_1^{1/2} - \log \gamma_2^{1/2} )  \varphi \, \dd x.
\end{gather*}
In particular, if we take any $ \varphi \in H^1_0(B) $ we see that
\begin{equation*}
\duality{q_1 - q_2}{\gamma_1^{1/2} \gamma_2^{1/2} \varphi} = - \int_{\R^n} \gamma_1^{1/2} \gamma_2^{1/2} \nabla ( \log \gamma_1^{1/2} - \log \gamma_2^{1/2} )\cdot \nabla   \varphi \, \dd x,
\end{equation*}
which means that $ \log \gamma_1^{1/2} - \log \gamma_2^{1/2} \in H^1_0(B) $ is a weak solution of
\[ - \nabla \cdot \left( \gamma_1^{1/2} \gamma_2^{1/2} \nabla ( \log \gamma_1^{1/2} - \log \gamma_2^{1/2} ) \right) = \gamma_1^{1/2} \gamma_2^{1/2} (q_2 - q_1) \]
in $ B $. The well-posedness of the above elliptic divergence-form equation implies
\begin{equation}
\norm{\log \gamma_1 - \log \gamma_2}{}{H^1(B)} \lesssim \norm{q_1 - q_2}{}{H^{-1}(B)} \lesssim \norm{q_1 - q_2}{}{H^{-1}(\R^n)}. \label{es:FORWARDstability} \footnote{From now on, the implicit constants also depend on the a priori bound $ M $.}
\end{equation}
So the goal now is to bound the right hand side of last inequality by the boundary data. In order to do so, we are going to use the equivalent norm in $ H^{-1}(\R^n) $ given in terms of the Fourier transform --from now on, this will be denoted by $ \mathcal{F} $. Let $ t \geq 1 $ be a constant to be chosen later. One has
\begin{gather*}
\norm{q_1 - q_2}{2}{H^{-1}(\R^n)} \lesssim \int_{ \{ |k| < t\} } (1+|k|^2)^{-1} | \mathcal{F}(q_1-q_2) (k)|^2 \,\dd k \\
+ \int_{ \{ |k| \geq t\} } (1+|k|^2)^{-1} | \mathcal{F}(\nabla \log( \gamma_1 \gamma_2 ) \cdot \nabla ( \log \gamma_1 - \log \gamma_2 )) (k)|^2 \,\dd k \\
+ \int_{ \{ |k| \geq t\} } (1+|k|^2)^{-1} | \mathcal{F}(\Delta( \log \gamma_1 - \log \gamma_2 )) (k)|^2 \,\dd k.
\end{gather*}
Here we have made a distinction between low and high frequencies in order to take advantage of the decay of $ (1+|k|^2)^{-1} $ for $ |k| > t $, as $ t $ becomes large. This can be obviously done in the second term of the right hand side with the available smoothness for $ \log \gamma_j $. However, this does not seem to be the case for the last one. To avoid this we will again use the approximation to the identity introduced in Section \ref{sec:CGOreview}. Hence
\begin{gather*}
\int_{ \{ |k| \geq t\} } (1+|k|^2)^{-1} | \mathcal{F}(\Delta( \log \gamma_1 - \log \gamma_2 )) (k)|^2 \,\dd k \lesssim \\
\lesssim \int_{ \{ |k| \geq t\} } (1+|k|^2)^{-1} | \mathcal{F}(\nabla \cdot ( \psi_h \ast \nabla( \log \gamma_1 - \log \gamma_2 ))) (k)|^2 \,\dd k \\
+ \int_{ \{ |k| \geq t\} } (1+|k|^2)^{-1} |k|^2 | \mathcal{F}(\nabla( \log \gamma_1 - \log \gamma_2 ) - \psi_h \ast \nabla( \log \gamma_1 - \log \gamma_2 )) (k)|^2 \,\dd k \\
\lesssim \frac{1}{h^2 t^2} + \sum_{j = 1}^n \norm{\partial_{x_j}( \log \gamma_1 - \log \gamma_2 ) - \psi_h \ast \partial_{x_j}( \log \gamma_1 - \log \gamma_2 )}{2}{L^2(\R^n)}\\
\lesssim h^{-2} t^{-2} + h^{2 \varepsilon}.
\end{gather*}
In the last inequality we used that $ \gamma_j \in C^{1, \varepsilon} (\R^n) $ , $ \supp (\gamma_j - 1) \subset \overline{B} $ and $ \gamma_j (x) \geq 1/(2M) $ for all $ x \in \R^n $. Choosing now $ h = t^{-1/(1 + \varepsilon)} $ one sees that
\begin{equation*}
\int_{ \{ |k| \geq t\} } (1+|k|^2)^{-1} | \mathcal{F}(\Delta( \log \gamma_1 - \log \gamma_2 )) (k)|^2 \,\dd k \lesssim t^{-2\varepsilon / (1 + \varepsilon)}.
\end{equation*}
Hence,
\begin{gather}
\norm{q_1 - q_2}{}{H^{-1}(\R^n)} \lesssim t^{n/2} \sup_{|k| < t} | \mathcal{F}(q_1-q_2) (k) | + t^{-\varepsilon / (1 + \varepsilon)}. \label{es:phiq_1q_2}
\end{gather}

Our next step is to bound uniformly $ | \mathcal{F}(q_1-q_2) (k) | $ for $ |k| < t $. In order to achieve this, we are going to use the CGOs constructed by Haberman and Tataru. Indeed, we shall substitute into the estimate \eqref{es:IntegralboundC1} the solutions $u_{\zeta_j}=\gamma_j^{-1/2}e^{\rho_j}(1+r_{\zeta_j})$ with $j \in \{1,2\}$. These solutions exist for all $ |\zeta_j| \gtrsim 1 $, where the implicit constant depends on $ M $. Here $\rho_j(x)=\zeta_j\cdot x$ and
\[ \zeta_j := s \eta_j + i \left( \frac{k}{2} + \left( s^2 - \frac{|k|^2}{4} \right)^{1/2} \kappa_j \right), \]
with $ k $, $ \eta_j $ and $ \kappa_j $ as in (\ref{def:vectorZETA}) and the additional condition that $ \eta_1 = - \eta_2 $, $ \kappa_1 = - \kappa_2 $ and $ s \in [ \max(1,|k|/2), + \infty ) $. Thus, we have the following
\begin{gather}
\int_{\R^n} \nabla \gamma_2^{1/2} \cdot \nabla (\gamma_2^{-1/2} v_1 v_2) \, \dd x - \int_{\R^n} \nabla \gamma_1^{1/2} \cdot \nabla (\gamma_1^{-1/2} v_1 v_2) \, \dd x = \nonumber \\
= \duality{q_1 - q_2}{e^{i x \cdot k}} + \duality{q_1 - q_2}{e^{i x \cdot k}(r_{\zeta_1} + r_{\zeta_2})} + \duality{m_{q_1} r_{\zeta_1}}{e^{i x \cdot k} r_{\zeta_2}}  \label{eq:identitypot} \\
- \duality{m_{q_2} r_{\zeta_2}}{e^{i x \cdot k} r_{\zeta_1}}. \nonumber
\end{gather}
Let us proceed to bound the terms appearing on the right hand side of \eqref{eq:identitypot}. Let $ j $ and $ l $ belong to $ \{ 1, 2 \} $ and $ l \neq j $. On the one hand
\begin{gather}
\left| \duality{q_j}{e^{i x \cdot k}(r_{\zeta_j} + r_{\zeta_l})} \right| = \left| \duality{q_j}{\phi e^{i x \cdot k}(r_{\zeta_j} + r_{\zeta_l})} \right| \leq \nonumber \\
\leq \norm{q_j}{}{X^{-1/2}_{\zeta_j}} \norm{\phi e^{i x \cdot k} r_{\zeta_j}}{}{X^{1/2}_{\zeta_j}} + \norm{q_j}{}{X^{-1/2}_{\zeta_l}} \norm{\phi e^{i x \cdot k} r_{\zeta_l}}{}{X^{1/2}_{\zeta_l}} \nonumber \\
\lesssim \langle k \rangle^{1/2} \left( \norm{q_j}{}{\dot{X}^{-1/2}_{\zeta_j}} \norm{r_{\zeta_j}}{}{\dot{X}^{1/2}_{\zeta_j}} + \norm{q_j}{}{\dot{X}^{-1/2}_{\zeta_l}} \norm{r_{\zeta_l}}{}{\dot{X}^{1/2}_{\zeta_l}} \right). \label{es:REMINDERq}
\end{gather}
In the last inequality we used \eqref{es:multiBYphi}. On the other hand
\begin{equation}
\left| \duality{m_{q_j} r_{\zeta_j}}{e^{i x \cdot k} r_{\zeta_l}} \right| \lesssim \langle k \rangle^{1/2} \norm{r_{\zeta_j}}{}{\dot{X}^{1/2}_{\zeta_j}} \norm{r_{\zeta_l}}{}{\dot{X}^{1/2}_{\zeta_l}}. \label{es:bound2mq}
\end{equation}
Indeed, by (\ref{eq:Leibniz}) one has
\begin{gather*}
\left| \duality{m_{q_j} r_{\zeta_j}}{e^{i x \cdot k} r_{\zeta_l}} \right| = \left| \duality{m_{q_j} \phi r_{\zeta_j}}{\phi e^{i x \cdot k} r_{\zeta_l}} \right| \lesssim \\
\lesssim \norm{\phi r_{\zeta_j}}{}{L^2(\R^n)} \norm{\phi r_{\zeta_l}}{}{L^2(\R^n)} + \int_{\R^n} |\nabla(\phi r_{\zeta_j} \phi e^{i x \cdot k} r_{\zeta_l})| \,\dd x\\
\lesssim |\zeta_j|^{-1/2} \norm{\phi r_{\zeta_j}}{}{X_{\zeta_j}^{1/2}} |\zeta_l|^{-1/2} \norm{\phi r_{\zeta_l}}{}{X_{\zeta_l}^{1/2}} + \norm{\phi r_{\zeta_j}}{}{X_{\zeta_j}^{1/2}} \norm{\phi e^{i x \cdot k} r_{\zeta_l}}{}{X_{\zeta_l}^{1/2}},
\end{gather*}
where we followed the proof of Lemma 2.3 of \cite{HT} to bound the integral before the last inequality. Finally, (\ref{es:bound2mq}) follows from \eqref{es:multiBYphi}.

Now it is a consequence of \eqref{eq:identitypot}, \eqref{es:REMINDERq}, \eqref{es:bound2mq} and \eqref{es:IntegralboundC1} that
\begin{gather*}
\left| \duality{q_1 - q_2}{e^{i x \cdot k}} \right| \lesssim \left( \norm{\Lambda_{\gamma_1} - \Lambda_{\gamma_2}}{}{} + \norm{\Lambda_{\gamma_1} - \Lambda_{\gamma_2}}{\varepsilon/(1+\varepsilon)}{} \right) \norm{u_{\zeta_1}}{}{H^1(B)} \norm{u_{\zeta_2}}{}{H^1(B)} \\
+ \langle k \rangle^{1/2} \sum_{l,j = 1}^2 \norm{q_l}{}{\dot{X}^{-1/2}_{\zeta_j}} \norm{r_{\zeta_j}}{}{\dot{X}^{1/2}_{\zeta_j}} + \langle k \rangle^{1/2} \norm{r_{\zeta_1}}{}{\dot{X}^{1/2}_{\zeta_1}} \norm{r_{\zeta_2}}{}{\dot{X}^{1/2}_{\zeta_2}},
\end{gather*}
where the implicit constant already depends on $ \varepsilon $. Moreover,
\[ \norm{u_{\zeta_j}}{}{H^1(B)} \lesssim e^{\sqrt{2}Rs} s, \]
because of \eqref{eq:normaH1}, \eqref{es:rBOUNDBYq}, \eqref{eq:boundq} and \eqref{es:NONdecay} and since $ |\zeta_j| = \sqrt{2} s $. Hence,
\begin{gather*}
\left| \duality{q_1 - q_2}{e^{i x \cdot k}} \right| \lesssim \left( \norm{\Lambda_{\gamma_1} - \Lambda_{\gamma_2}}{}{} + \norm{\Lambda_{\gamma_1} - \Lambda_{\gamma_2}}{\varepsilon/(1+\varepsilon)}{} \right) e^{2\sqrt{2}Rs} s^2 \\
+ \langle k \rangle^{1/2} \sum_{l,j = 1}^2 \norm{q_l}{}{\dot{X}^{-1/2}_{\zeta_j}} \norm{q_j}{}{\dot{X}^{-1/2}_{\zeta_j}} + \langle k \rangle^{1/2} \norm{q_1}{}{\dot{X}^{-1/2}_{\zeta_1}} \norm{q_2}{}{\dot{X}^{-1/2}_{\zeta_2}}.
\end{gather*}
Here we just used \eqref{es:rBOUNDBYq}. We next take average of last estimate in $ (s, \eta) \in [\lambda, 2 \lambda] \times S $ (here $ \lambda $ can not be smaller than $ \max(1,|k| / 2) $) and, using H\"older's inequality, we show that
\begin{gather*}
\left| \duality{q_1 - q_2}{e^{i x \cdot k}} \right| \lesssim \left( \norm{\Lambda_{\gamma_1} - \Lambda_{\gamma_2}}{}{} + \norm{\Lambda_{\gamma_1} - \Lambda_{\gamma_2}}{\varepsilon/(1+\varepsilon)}{} \right) e^{4\sqrt{2}R\lambda} \lambda^2 \\
+ \langle k \rangle^{1/2} \sum_{m,j = 1}^2 \left( \frac{1}{\lambda} \int_S \int_\lambda^{2\lambda} \norm{q_m}{2}{\dot{X}^{-1/2}_{\zeta_j}} \, \dd s \, \dd l \right)^{1/2} \left( \frac{1}{\lambda} \int_S \int_\lambda^{2\lambda} \norm{q_j}{2}{\dot{X}^{-1/2}_{\zeta_j}} \, \dd s \, \dd l \right)^{1/2} \\
+ \langle k \rangle^{1/2} \left( \frac{1}{\lambda} \int_S \int_\lambda^{2\lambda} \norm{q_1}{2}{\dot{X}^{-1/2}_{\zeta_1}} \, \dd s \, \dd l \right)^{1/2} \left( \frac{1}{\lambda} \int_S \int_\lambda^{2\lambda} \norm{q_2}{2}{\dot{X}^{-1/2}_{\zeta_2}} \, \dd s \, \dd l \right)^{1/2}.
\end{gather*}
Since $ \gamma_j \in C^{1, \varepsilon} (\R^n) $ , $ \supp (\gamma_j - 1) \subset \overline{B} $ and $ \gamma_j (x) \geq 1/(2M) $ for all $ x \in \R^n $, we choose $ h = \lambda^{-1/(2 + 2\varepsilon)} $ in \eqref{eq:potentiallarge} and we see that, for $ \lambda \geq \max(1,|k|) $,
\begin{gather*}
\left| \duality{q_1 - q_2}{e^{i x \cdot k}} \right| \lesssim \left( \norm{\Lambda_{\gamma_1} - \Lambda_{\gamma_2}}{}{} + \norm{\Lambda_{\gamma_1} - \Lambda_{\gamma_2}}{\varepsilon/(1+\varepsilon)}{} \right) e^{4\sqrt{2}R\lambda} \lambda^2 \\
+ \langle k \rangle^{1/2} \left( \lambda^{-1} + \lambda^{-\varepsilon/(1 + \varepsilon)} + (1 + \langle k \rangle^2 / \lambda) \lambda^{-\varepsilon/(1 + \varepsilon)} \right).
\end{gather*}
In order to make a uniform bound for $ |k| < t $, we only have to consider $ \lambda \geq t\geq 1 $ and check that
\begin{gather*}
\sup_{|k| < t} \left| \mathcal{F}(q_1-q_2) (k) \right| \lesssim \left( \norm{\Lambda_{\gamma_1} - \Lambda_{\gamma_2}}{}{} + \norm{\Lambda_{\gamma_1} - \Lambda_{\gamma_2}}{\varepsilon/(1+\varepsilon)}{} \right) e^{4\sqrt{2}R\lambda} \lambda^2\\
+ t^{3/2} \lambda^{-\varepsilon/(1 + \varepsilon)}.
\end{gather*}
This is the uniform bound of $ | \mathcal{F}(q_1-q_2) (k) | $ with $ |k| < t $ that we were looking for.

We now go back to \eqref{es:phiq_1q_2}, we plug in the uniform bound of $ | \mathcal{F}(q_1-q_2) (k) | $ with $ |k| < t $ and we choose $ \lambda = t^{1 + (n / 2 + 3 / 2)(1 + \varepsilon) / \varepsilon} $ (note that this choice of $ \lambda $ still satisfies the condition $ \lambda \geq t $ for $ t \geq 1 $). Since $ (n / 2 + 3 / 2)(1 + \varepsilon) < n + 3 \leq 2n $, we have
\begin{gather*}
\norm{q_1 - q_2}{}{H^{-1}(\R^n)} \lesssim  \left( \norm{\Lambda_{\gamma_1} - \Lambda_{\gamma_2}}{}{} + \norm{\Lambda_{\gamma_1} - \Lambda_{\gamma_2}}{\varepsilon/(1+\varepsilon)}{} \right) e^{4\sqrt{2}Rt^{1 + 2n/\varepsilon}} t^{n/2} t^{2(1 + 2n/\varepsilon)}\\
+ t^{-\varepsilon/(1 + \varepsilon)}.
\end{gather*}
If we now choose
\[ t = \left( \frac{1}{2c} \log \norm{\Lambda_{\gamma_1} - \Lambda_{\gamma_2}}{-1}{} \right)^{\varepsilon / (\varepsilon + 2n)} \]
with $ c > 4\sqrt{2}R (1 + \varepsilon) / \varepsilon $, we get
\begin{equation}
\norm{q_1 - q_2}{}{H^{-1}(\R^n)} \lesssim \left( \log \norm{\Lambda_{\gamma_1} - \Lambda_{\gamma_2}}{-1}{} \right)^{- \varepsilon^2 / (5n)}. \label{es:H^-1stability}
\end{equation}
Note that this choice of $ t $ is only possible if
\[ \norm{\Lambda_{\gamma_1} - \Lambda_{\gamma_2}}{}{} \leq e^{-2c} \]
since $ t $ can not be smaller than $ 1 $. However, if $ \norm{\Lambda_{\gamma_1} - \Lambda_{\gamma_2}}{}{} > e^{-2c} $ the estimate \eqref{es:H^-1stability} trivially holds. It is now an immediate consequence of \eqref{es:H^-1stability} and \eqref{es:FORWARDstability} the first stability estimate
\begin{equation}
\norm{\gamma_1 - \gamma_2}{}{H^1(\R^n)} \lesssim \left( \log \norm{\Lambda_{\gamma_1} - \Lambda_{\gamma_2}}{-1}{} \right)^{- \varepsilon^2 / (5n)}. \label{es:H^1stability}
\end{equation}
Using that, for any $ f \in L^\infty (\R^n) $,
\[ |f(x)|^{n / (1 - \delta)} \leq \norm{f}{n / (1 - \delta) - 2}{L^\infty (\R^n)} |f(x)|^2 \]
almost every $ x \in \R^n $, one gets
\[ \norm{\gamma_1 - \gamma_2}{}{W^{1, n / (1 - \delta)}(\R^n)} \lesssim \norm{\gamma_1 - \gamma_2}{ (1 - \delta) 2/n}{H^1(\R^n)}. \]
Finally, the Morrey embedding allows to bound
\begin{equation}
\norm{\gamma_1 - \gamma_2}{}{C^{0,\delta} (\R^n)} \lesssim \norm{\gamma_1 - \gamma_2}{}{W^{1, n / (1 - \delta)}(\R^n)}. \label{es:Morrey's}
\end{equation}
Theorem \ref{th:stability} follows now from \eqref{es:H^1stability} and \eqref{es:Morrey's}.

\section{Final discussion} \label{sec:finalDISCUSSION}
We end these notes with a final discussion about the possibility of improving Theorem \ref{th:stability}. As a result of this, we pose two naive questions.

Haberman and Tataru's improvements of the classical method based on the construction of CGOs, allowed them to prove uniqueness of the Calder\'on problem for continuously differentiable conductivities. However, we think these improvements (as they appear in \cite{HT}) do not provide stability for the class of conductivities $ \gamma \in C^1(\overline{\Omega}) $ (the space of bounded continuous functions in $ \Omega $ such that their partial derivatives $ \partial^\alpha \gamma $ are bounded and uniformly continuous in $ \Omega $) satisfying $ \gamma (x) > 1/M $ for all $ x \in \Omega $ and $ \norm{\gamma}{}{C^1(\overline{\Omega})} < M $, for an a priori given constant $ M > 1 $. We think so for two different reasons. The first reason comes up when looking at \eqref{es:mqCHOOSEepsilon} since, in order to make the right hand side of \eqref{es:mqCHOOSEepsilon} small for any $ |\zeta| $ larger than some constant independent of $ \gamma $ in the class where one wants to prove stability, this class has to enjoy a property of equicontinuity for their derivatives. Note that this smallness of the right hand side of \eqref{es:mqCHOOSEepsilon} is necessary to construct the CGOs. The second reason is related to the remainder properties of the CGOs. When looking at \eqref{eq:potentiallarge}, one is forced to require the class, where one wants to get stability, to consist of functions whose partial derivatives have the same $ L^2 $-modulus of continuity. The importance of this is due to the following two facts: the remainder properties of the CGOs associated to a conductivity $ \gamma $ will be shaped by the $ L^2 $-modulus of continuity of $ \partial^\alpha \gamma $ with $ |\alpha| = 1 $, and the behaviours of the remainders have to be the same to get a stability estimate for a whole class. On the other hand, the remainder properties of the CGOs condition the modulus of continuity of the resulting stability. So, if one wants to reach the optimal $ \log $-type stability for the Calder\'on problem one has to assume a H\"older $ L^2 $-modulus of continuity for $ \partial^\alpha \gamma $ with $ |\alpha| = 1 $.

The assumptions in Theorem \ref{th:stability} are slightly stronger than the requirements described above because of two reasons. The first one is that the only boundary stability results for the gradients of the conductivities that we know are by Alessandrini, where the coefficients are assumed to be in $ C^{1, \varepsilon} (\overline{\Omega}) $ with $ \Omega $ having Lipschitz boundary (see \cite{A90}); and by Sylvester and Uhlmann, where the coefficients are assumed to be continuously differentiable in a domain having smooth boundary (see \cite{SU88}). The second reason is that the requirements pointed out above are not only related to the coefficient defined on $ \Omega $ but also to their extensions to $ \R^n $. Thus, in order to provide a result under these mere requirements, one has to carry out extensions of the conductivities that keep the equicontinuity, the H\"older $ L^2 $-modulus of continuity and such that the values of their extensions outside $ \Omega $ only depend on the values of the conductivities on $ \partial \Omega $. The extensions we performed do keep the equicontinuity and satisfy the last condition (see \eqref{es:controlOUTSIDE}), however, it does not preserve $ L^2 $-moduli of continuity.

The issue concerning the stability on the boundary might find an answer in a method due to Brown (see \cite{B01}). In his lectures in The Special Trimester on Inverse Problems in Madrid (2011), Brown showed that his method for recovering continuous conductivities on Lipschitz boundaries could be extended to recover the gradient of continuously differentiable conductivities also on Lipschitz boundaries.

We end these notes posing two naive questions motivated by this discussion:
\begin{enumerate}
\item Is it possible to find a class of admissible conductivities defined on $ \overline{\Omega} $ and perform extensions of that conductivities to $ \R^n $, from their values on $ \partial \Omega $, such that the partial derivatives of those extensions form an equicontinuous class and they belong to $ \Lambda^{2,\infty}_\varepsilon (\R^n) $ with $ 0 < \varepsilon < 1 $ (see \cite{S} for the definition of this space)?

\item Can Brown's method be extended to prove stability for the gradient of continuously differentiable conductivities on Lipschitz boundaries?
\end{enumerate}
We think that positive answers to these questions lead to a slight improvement of Theorem \ref{th:stability}.

\end{document}